\documentclass[12pt,a4paper]{amsart}
\usepackage{amsfonts}
\usepackage{amssymb}
\usepackage{bbm}
\usepackage{comment}
\usepackage{enumitem}
\usepackage{dsfont}
\usepackage{amsthm}
\usepackage{amsmath}
\usepackage{commath}
\usepackage{url}
\usepackage{epsfig}
\usepackage{subfigure}
\usepackage{adjustbox}

\makeatletter
\@namedef{subjclassname@2010}{%
	\textup{2010} Mathematics Subject Classification}
\makeatother

\usepackage[T1]{fontenc}

\newtheorem{theorem}{Theorem}

\newtheorem{rmk}{Remark}

\theoremstyle{definition}

\newcommand{\ls}{\leqslant}
\newcommand{\gs}{\geqslant}

\newcommand{\one}{\mathds{1}}
\newcommand{\xn}{(x_n)_{n \in \mathbb{N}} }
\newcommand{\yn}{(y_n)_{n \in \mathbb{N}} }
\newcommand{\zn}{(z_n)_{n \in \mathbb{N}} }

\def \bE {\mathbb E}
\def \bN {\mathbb N}
\def \bP {\mathbb P}

\def \bR {\mathbb R}
\def \bZ {\mathbb Z}
\def \cX {\mathcal{X}}
\def \cY {\mathcal{Y}}
\def \cA {\mathcal{A}}
\def \dn {\mathcal{D}_N}
\def \pmn {(m,n)}
\def \pkl {(k,\ell)}

\def \dd {\,\mathrm{d}}

\begin{document}
	\baselineskip=17pt
	\title{On Inhomogeneous Poissonian Pair Correlations}
	\author{Manuel Hauke}
	\address{University of York, UK}
	\email{manuel.hauke@york.ac.uk}
	\author[A. Zafeiropoulos]{Agamemnon Zafeiropoulos}
	\address{Technion, Israel Institute of Technology}
	\email{agamemn@campus.technion.ac.il}
	
	\date{}
	
	\begin{abstract} We study the notion of inhomogeneous Poissonian pair correlations, 
		proving several properties that show similarities and differences to its homogeneous counterpart. In particular, we show that sequences with inhomogeneous Poissonian pair correlations need not be uniformly distributed, contrary to what was till recently believed.
	\end{abstract}
	
	\maketitle

	\section{Introduction}
	Let $\xn \subseteq [0,1]$ be a sequence. The {\it pair correlation counting function} of $\xn$ is defined by \[ R_2(s,N) = \frac{1}{N}\#\Big\{ 1 \ls m \neq n \ls N: \lVert x_m -x_n \rVert \ls \frac{s}{N}  \Big\}, \qquad s>0,  \]
	where we write $\|x\| = \min\{ |x-k| : k\in \mathbb{Z} \}$ for the distance of $x\in \mathbb{R}$ to its nearest integer. The sequence $\xn$ is said to have {\it Poissonian pair correlations} (or PPC, for abbreviation) if $$\lim_{N\to \infty}R_2(s,N) = 2s \qquad  \text{ for all } s>0.$$ The notion of PPC is known to have connections with mathematical physics, which go far beyond the scope of this paper. As instance, we only mention  the famous Berry--Tabor Conjecture  \cite{berry_tabor}. PPC have recently attracted increasing attention from the purely theoretical point of view. In particular, it has been shown (partially independently) by various authors that sequences with PPC are necessarily uniformly distributed \cite{alp,  sigrid, marklof, steinerberger}. \par The notion of PPC has an ``inhomogeneous'' variant, which nonetheless is not equally well-studied.  Given $\gamma \in \bR,$ we say that a sequence $\xn$ has $\gamma$\,-\,PPC if
	\[\lim_{N \to \infty} R_2(\gamma;s,N) = 2s \qquad \text{ for all }s > 0, \]
	where we define \[R_2(\gamma;s,N) = \frac{1}{N}\#\Big\{ 1\ls m \neq n \ls N: \lVert x_m -x_n - \gamma\rVert \ls \frac{s}{N} \Big\}. \]
	Clearly, for $\gamma=0$ and more generally for any integer value of $\gamma $ in the above definition, one recovers the classical property of PPC. Since the property of $\gamma$\,-\,PPC is invariant under integer translations of $\gamma,$
	we will only consider values of $\gamma \in [0,1].$  When $0< \gamma < 1$, we may refer to the property of $\gamma$\,-\,PPC as {\it inhomogeneous Poissonian pair correlations with respect to $\gamma$}. Also since $\gamma$\,-\,PPC is equivalent to $(1-\gamma)$\,-\,PPC, it will be sufficient to restrict our attention to values $0< \gamma \ls \frac12.$  \par A first treatment of the notion of $\gamma$\,-\,PPC can be found in \cite{ramirez}. There, it is explained that like PPC, the property of $\gamma$\,-\,PPC with $\gamma\neq 0$ is a {\it pseudorandomness} property, in the sense that if $(X_n)_{n=1}^{\infty}$ is a sequence of i.i.d. random variables, all following the uniform distribution in $[0,1],$ then with probability $1$ the sequence $(X_n)_{n=1}^{\infty}$ has $\gamma$\,-\,PPC. In \cite{ramirez} it is also stated among other results that just like their homogeneous counterpart, $\gamma$\,-\,PPC is a property that is stronger than uniform distribution; in other words,  any sequence with $\gamma$\,-\,PPC is automatically uniformly distributed. \par In the first result of our paper, we prove that this statement is actually false! 
	
	\begin{theorem} \label{the_theorem}
		Let $0<\gamma <1$. Then there exists a sequence $\xn \subseteq [0,1]$  that has $\gamma$\,-\,PPC but is not uniformly distributed.
	\end{theorem}
	In addition to Theorem \ref{the_theorem}, we try to shed some light on the reason why PPC imply uniform distribution while $\gamma$\,-\,PPC with $\gamma\neq 0$ do not.  Among the several proofs of the fact that PPC imply uniform distribution that were mentioned earlier, Aistleitner, Lachmann and Pausinger in \cite{alp} prove a stronger statement that connects the limiting behaviour of the pair correlation function $R_2(s,N)$ with the asymptotic distribution function of the sequence $\xn.$  We say that the function $G:[0,1] \to \mathbb{R}$ is the {\it asymptotic distribution function} of the sequence $\xn$ if 
	\[ G(x) = \lim_{N\to \infty} \frac{1}{N}\#\{ n\ls N : 0 \ls x_n \ls x\} \qquad \text{ for all } 0 \ls x \ls 1. \] 
	
	\noindent Their result is the following. \newline
	
	\noindent {\bf Theorem (ALP): } {\it Let $\xn\subseteq [0,1]$ be a sequence with  asymptotic distribution function $G:[0,1] \to \mathbb{R}.$ Let \mbox{$F:[0,\infty)\rightarrow [0,\infty]$} be defined by 
		\[ F(s) = \liminf_{N \to \infty}\frac{1}{N} \# \Big\{1 \ls m\neq n \ls N : \|x_m - x_n \| \ls \frac{s}{N}  \Big\}, \qquad s>0. \]
		Then the following hold: 
		\begin{itemize}
			\item[(i)] If $G$ is not absolutely continuous, then $F(s)=\infty$ for all $s>0$.
			\item[(ii)] If $G$ is absolutely continuous and $g$ is the corresponding density function (that is, $g = G'$ almost everywhere), then 
			\begin{equation}\label{thm_b}  \limsup_{s \to \infty}\frac{F(s)}{2s} \gs \int_0^1 g(x)^2 \,\mathrm{d}x. \end{equation}
	\end{itemize} }
	\medskip
	This result\footnote{In the original result in \cite{alp}, the definition of $F$ has limit instead of liminf. Furthermore, the ALP Theorem can be straightforwardly adapted to allow the asymptotic function $G$ not to be unique, but we assume so to keep the notation simple.} indeed implies that PPC is a property stronger than uniform distribution. If $\xn$ is not u.d. mod $1,$ then the integral on the right-hand side of \eqref{thm_b} is strictly greater than $1$ and $\xn$ cannot have PPC.  \par

	For uniformly distributed sequences, the density function is $g(x)=1$ and thus the ALP Theorem implies that for any $\varepsilon > 0$, $\liminf_{N\to \infty}R_2(s,N) \gs 2s- \varepsilon$ for arbitrarily large values of the scale $s>0.$  In other words, among all uniformly distributed sequences, for those who have PPC, the quantity $$\limsup_{s \to \infty}\liminf_{N\to\infty}\frac{R_2(s,N)}{2s}$$ exhibits an extremal behaviour, in the sense that it has the minimal possible asymptotic size, and this extremal behaviour is not attainable for non-uniformly distributed sequences. \par
	Returning to the inhomogeneous PPC, the proof of Theorem \ref{the_theorem} makes use of the following statement, which establishes a connection between the density function $g$ and the limiting behaviour of $R_2(\gamma;s,N)$ in a probabilistic context.  We say that a random variable $X$ on some probability space is uniformly distributed with respect to the function $G:[0,1]\to \mathbb{R}$ if $\bP(X< t) = G(t)$ for all $t\in [0,1].$
	
	\begin{theorem}
		\label{wish_thm_}
		Let $\gamma\in\bR$ and $G:[0,1]\to \mathbb{R}$ be an absolutely continuous distribution function with corresponding density function $g\in L^2([0,1]).$ Let $(X_n)_{n\in\bN}$ be a sequence of independent random variables on some probability space $(\Omega, \Sigma, \bP)$ that are uniformly distributed with respect to $G$. Then, for the $\gamma$-pair correlation function of $(X_n)_{n\in\bN}$ we have almost surely,
		\begin{equation}\label{limit_}
			\lim_{N\to \infty} R_2(\gamma;s,N) = 2s \cdot \int_{0}^{1} g(x)g(x + \gamma)\,\mathrm{d}x \qquad \text{ for all } s>0.
		\end{equation}
	\end{theorem}
	\begin{rmk}
		In the rest of the paper, when $g:[0,1]\to \mathbb{R}$ is a density function as above and $\gamma\neq 0,$ we implicitly extend $g$ to the real numbers periodically mod $1$ and write $\int_0^1 g(x)g(x+\gamma)\dd x$ instead of $\int_0^1 g(x)g(\{x+\gamma\})\dd x$, which would be more accurate.
	\end{rmk}
	
	We note that for the specific choice $\gamma=0,$ Theorem \ref{wish_thm_} is the heuristic observation made in \cite[eqn. (2)]{alp}. Theorem \ref{the_theorem} follows straightforwardly from Theorem \ref{wish_thm_} once we find a non-constant density  $g$ such that $\int_{0}^{1} g(x)g(x + \gamma)\,\mathrm{d}x = 1$ for a fixed $\gamma$.  \par
	After establishing the connection between $R_2(\gamma;s,N)$ and the density function $g$ for random variables described in \eqref{limit_}, we suspected that an analogue of the ALP Theorem is also true for $\gamma\neq 0.$ That is, if $g$ is the distribution function of the sequence $\xn$ and we define 
	\begin{equation} \label{F_gamma} F_\gamma(s) = \liminf_{N \to \infty}\frac{1}{N} \# \Big\{ 1\ls m\neq n \ls N : \|x_m - x_n -\gamma \| \ls \frac{s}{N}  \Big\}, \qquad s>0,  \end{equation}
	then  \begin{equation} \label{ideally} \limsup_{s \to \infty}\frac{F_\gamma(s)}{2s} \gs \int_0^1 g(x)g(x+\gamma)\,\mathrm{d}x.  \end{equation} 
	However, it turns out that \eqref{ideally} is not true, either.
	
	\begin{theorem}\label{inhom_thmB}
		For any $0<\gamma <1$, there exists a sequence $\xn \subseteq [0,1]$ with asymptotic density function $g: [0,1] \to \mathbb{R}$ such that for the function $F_\gamma$ as in \eqref{F_gamma} we have
		\[  \limsup_{s\to \infty} \frac{F_\gamma(s)}{2s} <  \int_0^1 g(x)g(x+\gamma)\,\mathrm{d}x. \]
	\end{theorem}
	
	Finally, we examine the relation between the property of $\gamma$\,-\,PPC for different values of $\gamma.$ Observe that in view of Theorem \ref{the_theorem}, $\gamma$\,-\,PPC for $\gamma \neq 0$ does not imply PPC;  otherwise, every sequence with $\gamma$\,-\,PPC would need to be uniformly distributed. 
	As a last result of the paper, we show that this is not a phenomenon that distinguishes between PPC and inhomogeneous pair correlations. In particular, it follows that the classical PPC property is not stronger than its inhomogeneous counterparts. 
	
	\begin{theorem} \label{thm3}
		Let $\gamma_1,\gamma_2 \in [0,1/2]$ with $\gamma_1 \neq \gamma_2$. Then there exists a sequence $\xn \subseteq [0,1]$ that has $\gamma_1$\,-\,PPC but not $\gamma_2$\,-\,PPC.   
	\end{theorem}
	
	\noindent {\bf Why $\gamma$\,-\,PPC?} Before proceeding to the proofs of the theorems, we think it would be worth discussing what we view as the motivation behind the definition of $\gamma$\,-\,PPC. Beyond pure interest in the notion itself, this motivation arises in an open problem from the metric theory of (homogeneous) Poissonian pair correlations. Given an increasing sequence $\cA = (a_n)_{n\in \bN} \subseteq \bN$, a series of results \cite{all,rs,rz,rz2} shows that under certain assumptions on $\cA$, the sequence $(a_nx)_{n\in \bN} $ has PPC for almost all $x\in [0,1]$. There also exist results \cite{all,walker} that provide conditions on $\cA$ under which the sequence $(a_n x)_{n\in \bN}$ has PPC for almost no $x\in [0,1].$ However, it remains an open question to determine whether for any choice of  $\cA$, the sequence $(a_n x)_{n\in \bN}$ has PPC either for almost all or for almost no $x\in [0,1],$ thus establishing a zero-one law in the theory of metric Poissonian pair correlations. \par For $\cA$ fixed, writing $ X_{\cA} = \{ x\in [0,1] : (a_nx)_{n\in\bN}  \text{ has PPC}\},$ the aforementioned problem is equivalent to determining whether $\lambda(X_{\cA})=0$ or $1$ for any choice of $\cA.$ To answer this question, it would suffice to check whether the set $X_\cA$ is invariant under the ergodic transformation $T(x) = 2x \mod 1.$  Given $x\in X_\cA,$ we have $2x\in X_\cA$ if and only if the pair correlation function of the sequence $(2a_nx)_{n\in\bN}$ satisfies $\lim_{N\to \infty}R_2(2x,2s,N) = 4s$ for all $s>0$. But 
	\begin{align*} R_2(2x,2s,N) &= \frac{1}{N}\#\Big\{ 1\ls m\neq n \ls N : \|a_nx-a_mx\| \ls \frac{s}{N}\Big\} \\  & \qquad \quad +\frac{1}{N}\#\Big\{ 1\ls m\neq n \ls N : \|a_nx-a_mx-\tfrac12 \| \ls \frac{s}{N}\Big\}, \end{align*} and therefore $2x\in X_\cA$ if and only if the sequence $(a_nx)_{n\in\bN}$ has $\frac{1}{2}$-PPC! In particular, Theorem \ref{thm3} shows that the na\"{i}vest approach of immediately implying $2x \in X_A$ does not work and thus, if possible at all, one has to use the additional structure of the sequence.

	\medskip
	
	\noindent {\bf Notation.}
	Given $x\in \mathbb{R}$ and $r>0$ we write $B(x,r)=\{ t \in \mathbb{R} : \|t-x\| \ls r\}$ for the interval of points in the unit torus that have distance at most $r$ from $x.$  Also throughout the text we shall use the standard Vinogradov $\ll$-notation: we write $f(x)\ll g(x), x\to \infty$ when $\limsup_{x\to \infty} f(x)/g(x) < \infty$.  \\

	\noindent {\bf Acknowledgements.} MH is supported by the EPSRC grant EP/X030784/1. A part of this work was supported by the Swedish Research Council under grant no. 2016-06596 while MH was in residence at Institut Mittag-Leffler in Djursholm, Sweden in 2024. AZ is supported by European Research Council (ERC) under the European Union’s Horizon 2020 Research and Innovation Program, Grant agreement no. 754475.
	
	\section{Proof of Theorems \ref{the_theorem} and \ref{wish_thm_} }
	
	As explained in the introduction, we begin with the proof of Theorem \ref{wish_thm_}, which is a more general result of probabilistic nature. The existence of sequences $\xn$ with $\gamma$--PPC that are not uniformly distributed will then follow as a simple corollary.

	\subsection{Proof of Theorem \ref{wish_thm_}} When $m\neq n,$ the random variable $X_m - X_n$ has probability density function 
	\begin{align} \label{d_function}
		d(t)   &  = \int_0^1 g(x)g(x+t)\dd x.
	\end{align}
	
	The $\gamma$-pair correlation function 
	\[ R_2(\gamma; s, N) = \frac{1}{N}\sum_{m\neq n\ls N}\one_{B(0,\frac{s}{N})}(X_m-X_n-\gamma) = \frac{1}{N}\sum_{m\neq n\ls N} \one_{B(\gamma,\frac{s}{N})}(X_m-X_n)  \]
	is itself a random variable with expectation
	\begin{align*}
		\bE[R_2(\gamma;s,N)] &= \frac{1}{N}\sum_{m\neq n \ls N}\int \one_{B(\gamma,\frac{s}{N}) }(X_m-X_n)\dd \bP \\ & = \frac{1}{N} \sum_{m\neq n \ls N}\int_{B(\gamma,\frac{s}{N}) } d(t)\dd t = (N-1)\int_{B(\gamma,\frac{s}{N}) }d(t) \dd t.
	\end{align*}
	Since $g\in L^2([0,1]),$ the function $d(t)$ defined in \eqref{d_function} is bounded and continuous, hence 
	\[ \lim_{N\to\infty}N \int_{B(\gamma,\frac{s}{N}) } d(t)\dd t =  \lim_{N\to \infty}2s\cdot \frac{1}{2s/N} \int_{B(\gamma,\frac{s}{N}) }d(t)\dd t = 2s\cdot d(\gamma). \]
	This allows us to conclude that \[ \lim_{N\to \infty}\bE[R_2(\gamma;s,N)] = 2s \cdot \int_{0}^{1} g(x)g(x + \gamma) \,\mathrm{d}x.  \]
	In the rest of the proof, we will write $\dn = \{ (m,n)\in \bN^2 : 1\ls m \neq n \ls N\}$ for the set of pairs of indices appearing in the pair correlation functions. The second moment of $R_2(\gamma;s,N)$ is equal to 
	\begin{equation} \label{second_moment}
		\begin{split}   \bE[R_2(\gamma;s,N)^2] = \frac{1}{N^2}\sum_{\pmn \in \dn }\int \one_{B(0,\frac{s}{N})}(X_m-X_n-\gamma) \dd \bP \qquad \qquad \qquad \qquad \\
			\qquad \qquad  + \frac{1}{N^2}\mathop{\mathop{\sum\sum}_{\pmn, \pkl \in \dn}}_{(k,\ell) \neq (m,n)}\hspace{-2mm} \int \one_{B(0,\frac{s}{N}) }(X_k-X_\ell-\gamma)\one_{ B(0,\frac{s}{N})}(X_m-X_n-\gamma) \dd \bP. \end{split}
	\end{equation}
	The first of the terms above is equal to $\frac{1}{N} \bE[ R_2(\gamma;s,N)].$ In the second sum appearing in \eqref{second_moment}, the contribution of pairs $(k,\ell), (m,n) \in \dn$ that share one common coordinate is 
	\begin{multline*}
		\frac{1}{N^2}\mathop{\mathop{\sum\sum}_{k,\ell,m \ls N}}_{\text{distinct}} \int \one_{B(0,\frac{s}{N})}(X_k-X_\ell-\gamma) \one_{B(0,\frac{s}{N})}(X_k-X_m-\gamma) \dd \bP \quad + \\ +  \frac{1}{N^2}\mathop{\mathop{\sum\sum}_{k,\ell,m \ls N}}_{\text{distinct}} \int \one_{B(0,\frac{s}{N})}(X_k-X_\ell-\gamma) \one_{B(0,\frac{s}{N})}(X_m-X_k-\gamma) \dd \bP, 
	\end{multline*}
	which is equal to  
	\begin{multline*}
		\frac{1}{N^2}\mathop{\mathop{\sum\sum}_{k,\ell,m \ls N}}_{ \text{distinct}} \int \one_{B(X_k-\gamma,\frac{s}{N})}(X_\ell) \one_{B(X_k-\gamma,\frac{s}{N})}(X_m) \dd \bP \quad + \\ +  \frac{1}{N^2}\mathop{\mathop{\sum\sum}_{k,\ell,m \ls N}}_{\text{distinct}} \int \one_{B(X_k-\gamma,\frac{s}{N})}(X_\ell) \one_{B(X_k+\gamma,\frac{s}{N})}(X_m) \dd \bP.
	\end{multline*}
	Since for $\ell \neq m$ the random variables $X_\ell, X_m$ are independent, the above is equal to 
	\begin{multline*}
		\frac{1}{N^2}\mathop{\mathop{\sum\sum}_{k,\ell,m \ls N}}_{ \text{distinct}} \int \one_{B(X_k-\gamma,\frac{s}{N})}(X_\ell)\dd\bP \int \one_{B(X_k-\gamma,\frac{s}{N})}(X_m) \dd \bP \quad + \\ +  \frac{1}{N^2}\mathop{\mathop{\sum\sum}_{k,\ell,m \ls N}}_{\text{distinct}} \int \one_{B(X_k-\gamma,\frac{s}{N})}(X_\ell) \dd \bP \int \one_{B(X_k+\gamma,\frac{s}{N})}(X_m) \dd \bP = \\  =  \frac{1}{N^2}\mathop{\mathop{\sum\sum}_{k,\ell,m \ls N}}_{\text{distinct}} \int \one_{B(0,\frac{s}{N})}(X_k-X_\ell-\gamma) \dd\bP \int \one_{B(0,\frac{s}{N})}(X_k-X_m-\gamma) \dd \bP  + \\ +  \frac{1}{N^2}\mathop{\mathop{\sum\sum}_{k,\ell,m \ls N}}_{\text{distinct}} \int \one_{B(0,\frac{s}{N})}(X_k-X_\ell-\gamma)\dd\bP \int  \one_{B(0,\frac{s}{N})}(X_m-X_k-\gamma) \dd \bP.
	\end{multline*}
	Back to \eqref{second_moment}, for pairs $(k,\ell), (m,n) \in \dn$ with no common coordinate, the random variables $X_k-X_\ell$, $X_m-X_n$ are independent and the corresponding contribution is \begin{multline*}
		\frac{1}{N^2}\mathop{\mathop{\sum\sum}_{k,\ell,m,n\ls N}}_{k,\ell \notin \{m,n\}}\int  \one_{B(0,\frac{s}{N})}(X_k-X_\ell-\gamma)\one_{B(0,\frac{s}{N})}(X_m-X_n-\gamma)\dd\bP = \\ = \frac{1}{N^2}\mathop{\mathop{\sum\sum}_{k,\ell,m,n\ls N}}_{k,\ell \notin \{m,n\}}\int  \one_{B(0,\frac{s}{N})}(X_k-X_\ell-\gamma)\dd \bP \int \one_{B(0,\frac{s}{N})}(X_m-X_n-\gamma)\dd\bP.
	\end{multline*}
	If we apply the same case distinction to the sum in
	\begin{align*}
		\bE[R_2(\gamma;s,N)]^2 &= \frac{1}{N^2}\Big( \sum_{(m,n)\in \dn} \int \one_{B(0,\frac{s}{N})}(X_m-X_n-\gamma)\dd \bP \Big)^2
	\end{align*}
	and combine with the equations above, we deduce that 
	\begin{multline*}
		\int \Big| R_2(\gamma;s,N) - \bE[ R_2(\gamma;s,N)] \Big|^2 \dd \bP = \bE[ R_2(\gamma;s,N)^2 ] - \bE[ R_2(\gamma;s,N) ]^2 \\
		= \frac{1}{N}\bE[ R_2(\gamma;s,N)] - \frac{1}{N^2}\sum_{(m,n)\in \dn}\Big( \int \one_{B(\gamma,\frac{s}{N})}(X_m-X_n)\dd\bP\Big)^2 \\ \ls \frac{1}{N}\bE[ R_2(\gamma;s,N)]  \, \ll \,  \frac{1}{N}, \quad N\to \infty, 
	\end{multline*}
	where we used again that $g \in L^2([0,1])$.
	An application of Chebyshev's inequality in combination with the first Borel--Cantelli Lemma shows that almost surely we have 
	\[ \lim_{N\to \infty} R_2(\gamma;s,N^2) = \lim_{N\to \infty} \bE[ R_2(\gamma;s,N^2)]  = 2s \cdot \int_{0}^{1} g(x)g(x + \gamma) \,\mathrm{d}x,\] and a standard approximation argument shows that the same is true along the whole sequence $R_2(\gamma;s,N).$

	\subsection{ Proof of Theorem \ref{the_theorem}} In view of Theorem \ref{wish_thm_}, in order to prove the existence of a sequence with $\gamma$\,-\,PPC that is not uniformly distributed, it simply suffices to define a density function $g$ such that $\int_0^1 g(x)g(x+\gamma)\dd x =1$ and $g$ is not identically equal to $1.$ Then, any sequence of random variables $\xn$ that has probability density function equal to $g$ will almost surely have $\gamma$-PPC (by \eqref{limit_}) but will not be uniformly distributed; otherwise we would have $g(x)=1$ for all $x\in [0,1].$ \\
	
	Given $\gamma \neq 0$, we first consider the case when $0 < \gamma < 1/2$. Choose some $\delta> 0$ with
	$\delta < \gamma$ and $\delta < 1 - 2\gamma$ and define
	\[
	g(x) = \begin{cases}
		1/\sqrt{\delta} &\text{ if } x \in [0,\delta) \cup [\gamma,\gamma + \delta),\\
		\dfrac{1 - 2\sqrt{\delta}}{\gamma - \delta} &\text{ if } x \in [\delta,\gamma),\\
		0 &\text{ otherwise}.
	\end{cases}
	\]
	In the case when $\gamma = 1/2$, we take some $0 < \delta < 1/2$ and let 
	\[
	g(x) = \begin{cases}
		1/\sqrt{2\delta} &\text{ if } x \in [0,\delta) \cup [\gamma,\gamma + \delta),\\
		\dfrac{1 - \sqrt{2\delta}}{\gamma - \delta} &\text{ if } x \in [\delta,\gamma),\\
		0 &\text{ otherwise}.
	\end{cases}
	\]
	In both cases, it is straightforward to check by elementary computations that the following three statements hold:
	\[ \int_{0}^1 g(x) dx = 1, \quad \int_{0}^1 g(x)g(x+ \gamma) dx = 1 \quad \text{ and } \quad  g \not\equiv 1. \]
	This concludes the proof of Theorem \ref{the_theorem}.

	\section{Proof of Theorem \ref{inhom_thmB}}
	
	We first prove the result for $\gamma = 1/2$ and then we explain  how the proof can be generalised to values $0<\gamma < 1/2$. In both cases, we shall make use of the binary van der Corput sequence $(c_n)_{n\in \bN}$ that is defined as follows (see also \cite{kuipers}): writing $n-1 = a_m(n)2^m + \ldots + a_1(n)2 + a_0(n)$ for the binary expansion of the integer $n-1$, the $n$-th term of $(c_n)_{n\in\bN}$ is the number \[ c_n = \frac{a_0(n)}{2} + \frac{a_1(n)}{2^2} + \ldots + \frac{a_m(n)}{2^{m+1}} \cdot \] 
	
	\noindent {\bf When $\gamma = 1/2.$} We shall construct a sequence $ \xn \subseteq [0,1]$ that
	on the one hand is uniformly distributed, which implies $\int_{0}^1 g(x) g(x-1/2) \,\mathrm{d}x = 1$, but on the other hand 
	\[F_{\frac12}(s) := \lim_{N\to \infty}\frac{1}{N} \# \Big\{ m,n \ls N : \|x_m - x_n -\tfrac12 \| \ls \frac{s}{N}  \Big\} = 0 \qquad \text{for all } s>0. \]
	We define the auxiliary sequences $\yn,\zn$ as follows: we set $y_n = \frac{1}{2} c_n$ and  \[ z_n = \frac{1}{2} + y_n +\frac{1}{3\cdot 2^N}  \qquad \text{whenever } 2^{N-1} < n \ls 2^{N}. \]  
	By these definitions, it is obvious that $\yn \subseteq [0,1/2]$ and $\zn \subseteq [1/2,1]$. \par  The sequence $\xn$ is constructed recursively, with the $N$--th step involving the definition of the terms $x_n, \, 2(N-1)2^{N-1} < n \ls 2N2^N.$ 
	At step $N = 1$, we set
	\[x_1 = y_1, \quad x_2 = z_1,\quad  x_3 = y_2, \quad  x_4 = z_2.\]
	Assume that for some $N\gs 1$ we have defined all points $x_n, 1\ls n\ls 2N2^{N}.$  Set  \[ (x_{2N2^N + 1}, \ldots, x_{2(N+1)2^{N}}) = (y_1,z_1,y_2,z_2,\ldots,y_{2^{N}},z_{2^{N}}) \]  The $2(N+1)2^N$-tuple $(x_{2(N+1)2^N+1},\ldots, x_{2(N+1)2^{N+1}})$ is then defined by concatenating the $2^{N+1}$ tuple \[ (y_{2^{N}+1},z_{2^{N}+1},\ldots,y_{2^{N+1}},z_{2^{N+1}})\] $N+1$-times with itself.  
	In that way, for each $N\gs 1$ the terms $x_1,\ldots, x_{2N2^N}$ contain $N$ copies of the points $y_1,\ldots,y_{2^N}$ and $N$ copies of the points $z_1,\ldots, z_{2^N}.$   \par It is straightforward to check that the sequence $\xn$ is uniformly distributed, so it remains to prove that $F_{\frac12}(s)=0$. 
	For any $N\gs 1$, set \begin{align*}
		Y_N &= \Big\{\frac{j}{2^{N+2}}: 0 \ls j < 2^{N+1}\Big\} \quad \text{ and }\\ 
		Z_N &=   \Big\{\frac{1}{2} + \frac{j}{2^{N+2}} + \frac{1}{3\cdot 2^k}: 0 \ls j < 2^{N+1}, 0 \ls k \ls N+1 \Big\} \cap \big[\tfrac12,1\big].
	\end{align*}
	We claim that 
	\begin{equation}\label{min_distance}  \Big\lVert a - b - \frac{1}{2}\Big\rVert \gs \frac{1/12}{2^N} \qquad \text{ whenever } \quad a,b \in Y_N \cup Z_N. \end{equation}
	This is obvious for $a,b \in Y_N$. Since
	\[\min  Z_N  = \frac{1}{2} + \frac{1}{3\cdot 2^{N+1}}\qquad \text{ and } \qquad  \max  Z_N  \ls 1,\]
	the inequality in \eqref{min_distance} also holds when $a,b \in Z_N$. It remains to check the case when $a \in Y_N, b \in Z_N$ (or vice-versa); then 
	\begin{equation*}
		\Big\lVert a - b - \frac{1}{2}\Big\rVert
		\gs \min_{0 \ls k \ls N+1} \min_{n \in \bZ}  \Big\lVert \frac{n}{2^{N+2}} -  \frac{1}{3\cdot 2^k}\Big\rVert \gs \frac{1}{3\cdot 2^{N+2}} = \frac{1/12}{2^{N}},
	\end{equation*}
	since $3$ is coprime to $2^N$. This proves \eqref{min_distance}.
	
	Now given $M \in \bN$, let $N = N(M)\gs 1$ be defined by
	\begin{equation*}
		2N\cdot 2^N < M \ls2 (N+1) 2^{N+1}.\end{equation*} We then have the inclusion
	\[\begin{split}\{x_n: 1 \ls n \ls M\} &\subseteq  \{y_n: 1 \ls n \ls 2^{N+1} \} \cup  \{z_n: 1 \ls n \ls 2^{N+1} \} \subseteq Y_N \cup Z_N.
	\end{split}
	\]
	Thus for a fixed value of $s > 0$, whenever
	$ \#\left\{a,b \in Y_N \cup Z_N: \lVert a- b - 1/2 \rVert \ls\frac{s}{2N\cdot 2^N}\right\}=0,$ then also $\#\left\{n, m \ls M: \lVert x_n - x_m - 1/2 \rVert \ls\frac{s}{M} \right\}=0$. Hence for all $N \gs 6s$, \eqref{min_distance} implies that $\#\left\{n, m \ls M: \lVert x_n - x_m - 1/2 \rVert \ls\frac{s}{M} \right\}=0$ and we obtain
	$F_{1/2}(s) = 0$.
	
	\medskip
	
	\noindent {\bf When $0 < \gamma < 1/2$.} We set  $\varepsilon = 2^{-i}$ where $i\gs 1$ is large enough such that $0 < \varepsilon < \min\{\tfrac{1}{2}(\tfrac{1}{2}-\gamma),\gamma\}$. We then define $\yn$ and $\zn$ by setting 
	\[ y_n = \varepsilon c_n \qquad \text{and} \qquad z_n = \gamma + y_n + \frac{1}{3\cdot 2^N}\quad \text{whenever } 2^{N-1} < n \ls 2^N.\] The sequence $\xn$ is then defined precisely as in the case  $\gamma = 1/2$.\par
	One can use the same arguments as before to show that for any $s>0$ we have $F_\gamma(s)=0.$ On the other hand, the asymptotic density function of $\xn$ is  $$g(x) = \begin{cases} 1/(2\varepsilon), & \text{ if } x\in [0,\varepsilon]\cup   [\gamma,\gamma + \varepsilon] \\ 0, & \text{ otherwise} \end{cases} $$ and thus
	$ \int_{0}^1 g(x)g(x+\gamma) \mathrm{d}x = 1/(4\varepsilon) > 0, $
	which concludes the proof. 
	
	\begin{rmk}
		Note that the sequence defined above for $\gamma \neq 1/2$ is not uniformly distributed. However, one could adapt the construction by ``diluting'' the sequence with i.i.d. samples on $[\varepsilon,\gamma] \cup [\gamma + \varepsilon,1]$ in order to get a sequence that almost surely is uniformly distributed mod $1$ and  satisfies
		\[\limsup_{s\to \infty} \frac{F_\gamma(s)}{2s} <  \int_0^1 g(x)g(x+\gamma)\,\mathrm{d}x = 1.\]
		We leave the details to the interested reader.
	\end{rmk}
	
	\section{$\gamma_1$\,-\,PPC does not imply $\gamma_2$\,-\,PPC}
	
	We finish with the proof of Theorem \ref{thm3}. Let $\xn$ be a sequence of independent random variables, following the uniform distribution in $[0,1]$. Define
	$\yn$ by \[ y_{2n-1} = x_n \quad \text{and} \quad y_{2n} = x_n +\gamma_2 \qquad (n\gs 1).\] 
	In what follows, we will write $R_2^{\cX}$ and $R_2^{\cY}$ for the pair correlation functions of $\xn$ and $\yn$ respectively. Then  by definition of $\yn$,
	\[ R_2^{\cY}(\gamma_2;s,2N) =  \frac{1}{2N}\#\Big\{n \neq m \ls 2N: \lVert y_n - y_m - \gamma_2 \rVert \ls \frac{s}{2N}\Big\}
	\gs \frac{1}{2N} N = \frac{1}{2} \cdot
	\]
	Taking  $s < \frac{1}{4}$, the above shows that $\yn$ does not have $\gamma_2$\,-\,PPC. We claim that $\yn$ has still almost surely $\gamma_1$\,-\,PPC: we have
	\begin{equation*}
		\begin{split}R_2^{\cY}(\gamma_1;s,2N)  &= \frac{1}{2N}\#\Big\{n\neq m \ls 2N: \lVert y_n - y_m -\gamma_1 \rVert \ls \frac{s}{2N}\Big\}
			\\& = \frac{1}{2N}\#\Big\{n \neq m \ls N: \lVert x_n - x_m -\gamma_1 \rVert \ls \frac{s/2}{N}\Big\}
			\\&\qquad +\frac{1}{2N}\#\Big\{n\neq m \ls N: \lVert x_n + \gamma_2 - (x_m + \gamma_2) - \gamma_1 \rVert \ls \frac{s/2}{N}\Big\}
			\\& \qquad + \frac{1}{2N}\#\Big\{n,m \ls N: \lVert x_n + \gamma_2 - x_m -\gamma_1 \rVert \ls \frac{s/2}{N}\Big\}
			\\& \qquad + \frac{1}{2N}\#\Big\{n,m \ls N: \lVert x_n - (x_m + \gamma_2) -\gamma_1 \rVert \ls \frac{s/2}{N}\Big\}
			\\&= R_2^{\cX}(\gamma_1;\tfrac{s}{2},N) + \tfrac12 R_2^{\cX}(\gamma_1-\gamma_2; \tfrac{s}{2}, N) + \tfrac12 R_2^{\cX}(\gamma_1+\gamma_2; \tfrac{s}{2}, N)
			\\&\qquad + \frac{1}{2N}\#\Big\{n \ls N: \lVert x_n + \gamma_2 - x_n -\gamma_1 \rVert \ls \frac{s/2}{N}\Big\}
			\\& \qquad + \frac{1}{2N}\#\Big\{n \ls N: \lVert x_n - (x_n + \gamma_2) -\gamma_1 \rVert \ls \frac{s/2}{N}\Big\}.
		\end{split}
	\end{equation*}
	Since $\gamma_1 - \gamma_2, \gamma_1 + \gamma_2 \notin \mathbb{Z}$, the last two terms in the above equation vanish if $N$ is sufficiently large.
	By Theorem \ref{wish_thm_}, the sequence $\xn$ has $\gamma_1$\,-\,PPC and $(\gamma_1-\gamma_2)$\,-\,PPC and $(\gamma_1+\gamma_2)$\,-\,PPC almost surely, therefore with probability $1$ we have
	\[\lim_{N \to \infty} R_2^\cY(\gamma_1;s,2N) = 2s \qquad \text{ for all } s>0. \]
	A standard approximation argument now gives that $\lim_{N\to \infty}R_2^\cY(\gamma_1;s,N) =2s$ for all scales $s>0$, which means that $\yn$ has $\gamma_1$\,-\,PPC almost surely.

\end{document}